\documentclass[11pt]{article}
\usepackage{latexsym}
\usepackage{amsmath}
\usepackage{amssymb}
\usepackage{amsthm}
\usepackage[pagebackref, plainpages=false, pdfpagelabels]{hyperref}

\renewenvironment{proof}{\noindent {\it Proof. }}{\hfill $\Box$}
\newtheorem {nummer} {} [section]

\newcommand{\Char} [1] {\mbox{char $#1$}}

\newcommand{\erz} [1] {\mbox{$\langle #1 \rangle$}}

\renewcommand {\phi} {\varphi}

\newcommand{\Rad} [1] {\mbox{{\rm Rad}$(#1)$}}
\newcommand{\Ref}[1]{\mbox{(}\ref{#1}\mbox{)}}

\newcommand {\N} {\mathbb{N}}

\newcommand {\K} {\mathbb{K}}
\newcommand {\F} {\mathbb{F}}

\begin{document}
\title{The geometry of hyperbolic lines in polar spaces}
\author{Hans Cuypers \\
Department of Mathematics \\
Eindhoven University of Technology \\
P.O. Box 513 \\
5600 MB Eindhoven \\
The Netherlands\\
f.g.m.t.cuypers@tue.nl}

\maketitle

\begin{abstract}
In this paper we consider partial linear spaces induced on the point set of a
polar space, but with as lines the hyperbolic lines of this polar space.
We give some geometric characterizations of these and related spaces.

The results have applications in group theory,
 in the theory of Lie algebras and in graph theory.
\end{abstract}

\section{Introduction}
Classical polar spaces arise from  symplectic, hermitian and  
pseudo-quadratic forms on a vector space $V$.
They consist of those 
points and lines of the 
corresponding projective space $\mathbb{P}(V)$ on which the form vanishes. 
Actually, the results of Tits \cite{tits}, Veldkamp \cite{veldkamp},
Buekenhout-Shult \cite{bu-sh} and Johnson \cite{johnson}
imply that 
under some weak additional conditions (on non-degeneracy and rank) 
all these  polar spaces can be characterized by the famous `one-or-all
axiom'; see \Ref{polardef}.
In this paper we  provide a characterization of the  
geometries that have the same point set as a polar
space, 
but whose lines are the so-called {\em hyperbolic} lines of the polar space.
Our results are inspired by and can be viewed as a geometric counterpart
of Timmesfeld's theory of transvection groups \cite{timmes,timmes-book}.
Before stating our results, we introduce some notation and terminology.

A {\em partial linear space} is a  
pair $(P,L)$ consisting of a set $P$ of {\em points}
and a set $L$ of subsets of $P$ of size at least $2$ called {\em lines}, such that
two points are in at most one line.
A partial linear space is called {\em thick} 
if all lines contain at least $3$ points.
If two points $p$ and $q$  are on a line, then they are called  {\em
  collinear}.
This is denoted by $p\sim q$, while the unique line containing  them
is denoted by $pq$.
By $p\perp q$ we denote that the points $p$ and $q$ are equal or
non-collinear.
(As we are interested in the geometry of hyperbolic lines, this notation
differs  from the common notation that $\perp$ denotes collinearity.)
By definition, a point $p$ is not collinear to itself. The set $p^\perp$,
respectively $p^\sim$, denotes the set of all points non-collinear, respectively
collinear to $p$. 
By $(P,\perp)$ and $(P,\sim)$ we denote the graphs with point set $P$ and 
two points $x,y$ being adjacent if and only if $x\perp y$ or $x\sim y$,
respectively. The space
$\Pi$ is called {\em connected} if $(P,\sim)$ is connected and {\em co-connected} if
$(P,\perp)$ is connected. 

A {\it subspace} $X$ of $\Pi$ is a subset
of $P$ with the property that any line intersecting 
it in at least two points is completely
contained in $X$. Since the intersection of subspaces is again a subspace,
we can, for each subset $Y$ of $P$, define the subspace $\langle Y \rangle $
{\it generated by } $Y$ to be the intersection of all subspaces containing $Y$.
A {\it plane} is a subspace generated by two intersecting lines.

A subspace $X$ is called {\it transversal} if 
for each line $\ell\subseteq X$ we have that $X$ is partitioned by the sets
$x^\perp\cap X$ for $x\in \ell$. The subsets $x^\perp \cap X$ are  called the 
{\it transversals} of $X$.
The following  result is the main theorem of this paper.

\begin{nummer}{\rm \bf Main Theorem}.\label{main}
Let $\Pi=(P,L)$ be a connected and co-connected  partial linear space satisfying
the following:
\begin{enumerate}
\item all lines in $L$ contain at least $4$ points;
\item if $x$ and $y$ are points with $x^\perp \subseteq y^\perp$ then
$x=y$;
\item every triple of points $x,y$ and $z$ with $x\sim y \sim z$ and $x\perp z$
generates a transversal plane;
\item for each point  $x\in P$ and transversal plane $\pi$ we have 
$x^\perp\cap\pi$ is empty, a point, a line, or  a transversal of $\pi$.
\end{enumerate}
Then $(P,\perp)$ is a non-degenerate polar graph of rank at least $2$.

Moreover, if
\begin{enumerate}
\item[{\rm (e)}] for all $h\in L$ we have $h=h^{\perp\perp}$,
\end{enumerate}
then $(P,L)$ is the geometry  of points and  hyperbolic lines of the
polar graph $(P,\perp)$.
\end{nummer}

In the geometry of points and hyperbolic lines of classical polar spaces
defined over fields (see \Ref{classicalplane}) 
any two intersecting hyperbolic 
lines are contained in a subspace
forming either a dual affine plane or a classical unital.
The following result characterizes the geometry of hyperbolic lines of these  
spaces by this property.
To state the result we need the following definition.
A partial linear space $\Pi=(P,L)$ is called {\it planar}, if and only if every pair of intersecting lines
is contained in a unique plane.

\begin{nummer}{\rm \bf Theorem}.\label{planethm}
Let $\Pi=(P,L)$ be a planar, connected and co-connected partial linear space.
Suppose the following hold in $\Pi$:

\begin{enumerate}
\item planes are either isomorphic to a dual affine plane or
are linear;
\item in a linear plane no $4$ lines intersect in $6$ points;
\item for all points $x$ and $y$ we have 
that $x^\perp \subseteq y^\perp$ implies $x=y$;
\item for all $l\in L$ we have $|l|\geq 4$ and $l^{\perp\perp}=l$.
\end{enumerate}
Then $\Pi$ is isomorphic to the geometry of hyperbolic lines in a non-degenerate
polar space of rank at least $2$.
\end{nummer}

The above theorem generalizes the characterizations of
the geometry of hyperbolic lines in symplectic spaces
as given by Hall \cite{copolar} and Cuypers  \cite{c-symp} to the 
geometry of hyperbolic lines in arbitrary polar spaces.

Notice that the classical unitals
do satisfy condition (b) of the above theorem.
This was already  observed by O'Nan for finite fields (see \cite{onan}), but also holds true in the infinite case; see \Ref{classicalplane}.
The restriction that lines contain at least $4$ points
is necessary, as the Fischer spaces of 3-transposition groups
with trivial $O_2$ 
do also satisfy the other conditions, see \cite{buek-Fischer,ch-class,ch-fischer,chps}.
In particular, the Fischer spaces of the sporadic groups 
${\rm Fi}_{22}$, ${\rm Fi}_{23}$ and ${\rm Fi}_{24}$  do satisfy
the conditions (a)-(c)  of the theorem  but have $3$ points per line.

Condition (b) is also necessary. 
Indeed, the generalized Fischer spaces of non-isotropic points
and tangent lines in
a unitary geometry  over $\F_4$ are geometries satisfying 
(a), (c) and (d), but not (b); see \cite{c-gen,cs}.
 
Combining the above results 
with \cite[Section 4]{chps} we obtain the following result.

\begin{nummer}{\rm \bf Theorem}.\label{fischer}
Let $\Pi=(P,L)$ be a  planar and connected partial linear space
of finite order $q\geq 2$.
Suppose the following hold in $\Pi$:

\begin{enumerate}
\item all planes are finite and either isomorphic to a dual affine plane or
a classical unital;
\item there exists a pair of  non-collinear points;
\item for all $x,y\in P$ we have that 
$x^\perp \subseteq y^\perp$ or $x^\sim\cup\{x\} \subseteq y^\sim\cup\{y\}$ implies $x=y$;
\item if $\pi$ is a linear plane of $\Pi$ and $x$ a point, then $x^\perp \cap \pi \neq \emptyset$.
\end{enumerate}
Then we have one of the following:

\begin{enumerate}
\item
$q$ is a prime power and $\Pi$ is isomorphic to the geometry of hyperbolic lines in a non-degenerate
symplectic or
unitary polar space over the field $\F_q$ or $\F_{q^2}$, respectively.
\item  $\Pi$ is a symmetric  Fischer space.
\item $\Pi$ is an orthogonal Fischer space over $\F_2$.
\item $\Pi$ is an orthogonal Fischer space over $\F_3$ of dimension at most 
$7$.
\item $\Pi$ is a  sporadic  Fischer space related to $\Omega^+(8,2):{\rm S}_3$ or 
${\rm Fi}_{22}$.
\end{enumerate}
\end{nummer}

The remainder of this paper is organized as follows.
In the next section we discuss various examples appearing in the 
theorems. Then in Section \ref{reduction-sec}
we start with the proof of Theorem \ref{main}.
Here we prove  Reduction Theorem
\ref{reduction} in which we show how a geometry satisfying
some weaker conditions than as in the hypothesis of Theorem \ref{main}
reduces to a geometry satisfying the hypothesis of \Ref{main}.
The proof of Theorem \ref{main} is then finished in Section \ref{geometry-sec}.
Section \ref{plane-sec} is   
devoted to a proof of Theorem \ref{planethm}, while Theorem
\ref{fischer} will be proven in Section \ref{finite-sec}.

Some of the results of this paper have been around for several
years and have been available in preprint-form.
Since over the years Theorem \ref{main} has 
found applications in group theory \cite{k-trans}, graph theory
\cite{altmann,gramlich} and in the theory of Lie algebras
\cite{ci,panhuis}, we have decided to publish the results.
In the final section
of this paper we briefly discuss these applications of our results.

\section{Hyperbolic lines in classical polar spaces}
\label{example-sec}

In this section we discuss some basic results on polar spaces
and their hyperbolic lines
and describe various  examples of geometries  occurring in the theorems
from the first section.

\begin{nummer}\label{polardef}
{\rm {\bf Definitions.}
A  {\em partial linear space}
 $\mathcal{S}=(\mathcal{P},\mathcal{L})$ 
is called a {\em polar space} if it satisfies the  so-called one-or-all axiom:

\begin{enumerate}
\item[] For each point $p$ and line $l$ of $\mathcal{S}$ with $p\not\in l$, 
the point $p$ is
collinear with either one or all points of $l$.
\end{enumerate}
The {\em radical} $\Rad{\mathcal{S}}$ of a polar space
$\mathcal{S}$ is the set of points which are collinear to all points.
A polar space is {\em non-degenerate}, if its radical is empty.
A non-degenerate polar space in which a point $p$ is either 
on a line $l$
or just collinear with a unique point of $l$  is also called 
a {\em generalized quadrangle}.

If $\mathcal{S}=(\mathcal{P},\mathcal{L})$ is a polar space and 
$p\in \mathcal{P}$, then by $p^\perp$ we denote the set of all points
equal or collinear to $p$. If $p,q\in \mathcal{P}$ then $p\perp q$
denotes that $q\in p^\perp$. (Notice the difference with the definition of
$\perp$ as given in the introduction!)
If $X$ is a subset of $\mathcal{P}$, then $X^\perp=\bigcap_{x\in X}\ x^\perp$.

The following two results are well known to experts in the field, but for 
the sake of completeness we include them in this paper.
}\end{nummer}

\begin{nummer}\label{non-degenerate}
If $\mathcal{S}=(\mathcal{P},\mathcal{L})$ is a non-degenerate polar space and
$x,y\in \mathcal{P}$ are points with $x^\perp\subseteq y^\perp$, then $x=y$. 
\end{nummer}

\begin{proof}
Suppose $x\neq y$ are points with  $x^\perp\subseteq y^\perp$.
Let $u$ be a point not collinear with $y$. Then $u$ is collinear
to a point $z$ on the line $xy$ through $x$ and $y$.
Let $v$ be a point not collinear to $z$. Then $v$ is collinear to a point 
$w$ on $xy$ different from $x,y$ and $z$.
The point $u$ is collinear to a point on $vw$, which we may assume to be $v$.
Now $y^\perp$ meets $uv$ in a point $r$. This point is then also
collinear with $x$ and hence with $z$ and $w$. 
But $u$ is the only point on $uv$ collinear with $z$ and $v$ the only point
collinear with $w$ and we find $u=r=v$ contradicting $u\neq v$.
\end{proof}

\medskip

A {\em hyperbolic line} of a polar space $\mathcal{S}$ is a 
set of points of the
form
$\{x,y\}^{\perp\perp}:=(\{x,y\}^\perp)^\perp$ where $x,y$ are two non-collinear points.

\begin{nummer}\label{hyperbolic}
Suppose that $\mathcal{S}$ is a thick non-degenerate polar space. 
Then  for any two points $z_1,z_2$ of a   hyperbolic line $H$
we have $z_1\not\perp z_2$ and $H=\{z_1,z_2\}^{\perp\perp}$.
\end{nummer}

\begin{proof}
Suppose $x,y\in \mathcal P$ are non-collinear points and $H=\{x,y\}^{\perp\perp}$
is the hyperbolic line on $x$ and $y$.
Let $z\in H$ be different from $x$ and $y$.
Then we have 
$x^\perp\cap y^\perp \subseteq z^\perp$.

First we prove that $x\not\perp z$.
Assume this is not true, so $z$ is collinear to $x$.
Every line $\ell$ on $x$ contains a point $x'\in y^\perp$ different from
$x$, which is then also in $z^\perp$. So the line $\ell$ contains two points
collinear to $z$ and hence  is contained in
$z^\perp$. 
But this implies that $x^\perp\subseteq z^\perp$, contradicting
$\mathcal S$ to be non-degenerate; see \Ref{non-degenerate}.
So, indeed $x\not\perp z$.

As $x^\perp\cap y^\perp \subseteq z^\perp$ we have $x^\perp\cap y^\perp\subseteq x^\perp \cap z^\perp$.
We will now prove these sets to be equal.
Suppose the contrary and let $u\in x^\perp\cap z^\perp$ but not in $y^\perp$.
Since $x$ and $z$ are not collinear,  $u$ is different from $x$.
The line on $x$ and $u$ contains a point $v\neq u$ in $y^\perp$.
But then $v\in x^\perp\cap y^\perp$ and hence also in $z^\perp$.
But now we find that the points $u$ and $v$ and hence also $x$, which is
on the line $uv$, are in $z^\perp$.
So indeed, $x^\perp\cap y^\perp= x^\perp \cap z^\perp$,
which implies that $H$ equals $\{x,z\}^{\perp\perp}$, which  is the hyperbolic line on $x$ and $z$.

Now suppose $z_1,z_2$ are two points in $H$, then applying the above twice
we find $H=\{x,z_1\}^{\perp\perp}=\{z_1,z_2\}^{\perp\perp}$ and $x\not\perp z_1\not\perp z_2$.
\end{proof}

\medskip

Denote by $\mathcal{H}$ the set of all hyperbolic lines
of a non-degenerate polar space $\mathcal{S}=(\mathcal{P},\mathcal{L})$.
Then $(\mathcal{P},\mathcal{H})$ is a partial linear space, called the
{\em geometry of hyperbolic lines} of $\mathcal{S}$.

In the remainder of this section we describe the classical polar spaces and 
their hyperbolic lines.

\begin{nummer} \label{classicalpolar}
{\rm {\bf Classical polar spaces.}
Let $\K$ be a skew field with anti-automorphism $\sigma$, 
where $\sigma^2=id$ and
$\epsilon \in \K$.
Assume that $W$ is a right vector space over $\K$.
As in  
Tits \cite[(8.2.1), (8.1.2)]{tits}, 
a $(\sigma,\epsilon)$-{\em hermitian form} $f$ on $W$ is a mapping
$f:W\times W \rightarrow \K$ such that for all $v,w,u\in W$
and $\lambda,\mu \in \K$ we have
\begin{enumerate}
\item $f(v+w,u)=f(v,u)+f(w,u)$ and $f(u,v+w)=f(u,v)+f(u,w)$ ;
\item $f(v\lambda,w\mu)=\lambda^\sigma f(v,w)\mu$;
\item $f(v,w)=f(w,v)^\sigma\epsilon$.
\end{enumerate}

The form $f$ is called {\em trace-valued} if and only if
$f(v,v)\in\{t+t^\sigma\epsilon\mid t\in \K\}$ for all $v\in V$.

We set
$$\K_{\sigma,\epsilon} := \{c - \epsilon c^\sigma \mid c\in \K \},\quad \K^{\sigma,\epsilon} := \{c \in \K \mid \epsilon c^\sigma = -c\}.$$

A {\em pseudo-quadratic form} (or $(\sigma,\epsilon)$-quadratic form)
$q$ with associated trace-valued $(\sigma,\epsilon)$-hermitian form
$f:W\times W\rightarrow \K$ is a map $q:W\rightarrow \K/\K_{\sigma,\epsilon}$, with the property
that
$$q(v\lambda)=\lambda^\sigma q(v)\lambda, {\ \rm and \ }q(v+w)= q(v)+q(w)+ f(v,w)+\K_{\sigma,\epsilon}$$
for all $v,w\in W$
and $\lambda\in \K$.

Suppose $W$ is endowed with a $(\sigma,\epsilon)$-hermitian form $f$
or pseudo-quadratic form $q$ with associated $(\sigma,\epsilon)$-hermitian form $f$.
For any subset $U$  of $W$, we set
$ U^\perp = \{w\in W \mid f(w,u)=0$ for all $u \in U\}$.
The {\em radical} of $f$ is the space $\Rad{W,f} = W^\perp$.
For any  $U,V\subseteq W$ 
we denote by $U\perp V$ that $f(u,v)=0$ for all $u\in U$ and $v\in V$.

A vector $w \in W$ with $f(w,w) = 0$ or
$q(w) = 0$, respectively, is called {\em isotropic}. (Note that such a
vector is also called {\em singular} in the literature.)
We  assume that there exists an isotropic vector, not contained in
$\Rad{W,f}$. A subspace $U$ of $W$ is called {\em isotropic} if $U\perp U$.
In the case where $W$ is endowed with a pseudo-quadratic form, we
define the {\em isotropic radical}
of $(W,q)$ to be the subspace   $R=\{r\in \Rad{W,f} \mid q(r)=0+\K_{\sigma,\epsilon}\}$.

We call $f$ {\em non-degenerate} if  $\Rad{W,f} =0$ 
and $q$ {\em non-degenerate}
if $R=\{0\}$.

Fix a complement $W_0$ to $R$ in $W$. Then $q|_{W_0}$ is non-degenerate.
We set
$$\Phi_q := \{c \in \K \mid \mbox{ there exists $r_c \in
\Rad{W_0,f}$ with $q(r_c) = c + \K_{\sigma,\epsilon}$} \}.$$ 
Then $\Phi_q$ is an additive subgroup of $\K$.
For each  $c \in \Phi_q$, there is a unique
$r_c\in \Rad{W_0,f}$ with $q(r_c)=c+\K_{\sigma,\epsilon}$.
Further, $d^\sigma \Phi_q d = \Phi_q$ for $0 \neq d \in \K$ and
$\K_{\sigma,\epsilon} \subseteq \Phi_q\subseteq \K^{\sigma,\epsilon}$.
If $\Rad{W_0,f} = 0$ (in particular if $\Char K \neq 2$), 
then $\Phi_q= \K_{\sigma,\epsilon}$.

By $\mathcal{S}=\mathcal{S}_f$ (or $\mathcal{S}_q$, respectively) we denote the polar space 
consisting 
of the point set $\mathcal{P}$ of isotropic $1$-spaces and line set
$\mathcal{L}$ 
of all isotropic $2$-spaces of $W$ with respect to $f$ 
(or $q$, respectively). 
Any  polar space arising in this way is called {\em classical}.

By the results of Tits \cite[(8.1) and (8.2)]{tits},
any non-degenerate classical polar space with non-empty line set
can be obtained (up to isomorphism) 
as $\mathcal{S}_f$ or $\mathcal{S}_q$ with $f$ or $q$ as below.

\begin{enumerate}
\item[(I)]  $f: W \times W \to \K$ is $(\sigma,-1)$-hermitian form
and $\K_{\sigma,-1}=\K^{\sigma,-1}$; 
\item[(II)] $q:W\rightarrow \K$ is a quadratic form, i.e. a $(id,1)$-quadratic
  form and $1\in \Phi_q$ or $\Phi_q=\{0\}$.
\item[(III)]
$q: W \to \K/\K_{\sigma,-1}$ is  
a pseudo-quadratic form  with associated (trace-valued)
$(\sigma,-1)$-hermitian form $f: W \times W \to \K$.
Moreover, $1\in \Phi_q$.
\end{enumerate}
}\end{nummer}

\begin{nummer} \label{hyperboliclines}
{\rm 
{\bf Hyperbolic lines in classical polar spaces.}
Let $\mathcal{S}$ be a non-degenerate classical
polar space $\mathcal{S}_f$ or $\mathcal{S}_q$
of some non-degenerate 
$(\sigma,\epsilon)$-hermitian form $f$ or (pseudo-)quadratic form $q$,
as in \Ref{classicalpolar}(I)-(III).  
Assume that the points of $\mathcal S$ generate $W$ and that $\mathcal S$
contains lines.

Let $v_1$ and $v_2$ be two isotropic vectors such that
$f(v_1,v_2) = 1$.
Then 
$x=\langle v_1\rangle$ and $y=\langle v_2\rangle$ are 
two non-collinear points of $\mathcal{S}$.   
The $2$-dimensional space $H := \erz{v_1,v_2}$ 
is non-degenerate and $W = H \perp H^\perp$.

\begin{nummer}
The hyperbolic line $\{x,y\}^{\perp\perp}$ equals
 $\{\langle v\rangle\mid v\in \langle v_1,v_2\rangle\}\cap \mathcal{S}_f$
or $\{\langle v\rangle\mid v\in \langle v_1,v_2\rangle+{\rm Rad}(W_0,f)\}\cap \mathcal{S}_q$.
\end{nummer}

\begin{proof}
Suppose $z\in \{x,y\}^{\perp\perp}$.
As  $\mathcal S$ is generated by $x^\perp$ and $y$, (and similarly $y^\perp$
and $x$), see \cite{cs}, we see that $\langle \{x,y\}^{\perp} \rangle =H^\perp$.
But that implies that $\{x,y\}^{\perp\perp}$ is contained in $H+{\rm Rad}(W_0,f)$.
\end{proof}

\medskip

Assume that $\mathcal{S}=\mathcal{S}_f$ 
for some $(\sigma,-1)$-hermitian
form $f$ as in (I). 
Then the isotropic $1$-spaces inside $H$ are $\langle v_2\rangle$
and  the spaces
$\langle v_1\lambda+  v_2  \rangle$ where $\lambda\in\K_{\sigma,-1}=\K^{\sigma,-1}$.
Indeed, 
$f( v_1\lambda+v_2, v_1\lambda+v_2)=
\lambda^\sigma-\lambda=0$ if and only if $\lambda^\sigma =\lambda$.

In particular,   as $1\in \K_{\sigma,-1}\neq \{0\}$, 
the hyperbolic lines are thick.
Denote by $\mathcal{H}$ the set of hyperbolic lines
of $\mathcal{S}$. Then $\Pi=(\mathcal{P},\mathcal{H})$ is 
a partial linear space.

Now suppose $v_3\in V\setminus H$ is an isotropic vector
with $f(v_3,v_2)=1$ and $f(v_1,v_3)=0$.   
Then $v_2+v_3$ is also an isotropic vector.
Now consider the non-degenerate $2$-spaces
$H_1=\langle v_1,v_2+v_3\rangle$ and $H_2=\langle  v_1\lambda +v_2,
 v_3\mu+v_2\rangle$,
where $\lambda\neq 0$ and $\mu\neq 0,1$ are distinct elements in $\K_{\sigma,-1}$.
The intersection $H_1\cap H_2$ is the $1$-space 
$\langle v_1\lambda(1-\mu^{-1})+v_2+v_3\rangle$.
This is a point of $\mathcal{S}$ if and only if
$\lambda\mu^{-1}\in \K_{\sigma,\epsilon}$.
So, if inside $\Pi$ 
 the three points $\langle v_1\rangle$, $\langle v_2\rangle$ and 
$\langle v_3\rangle$ are in a subspace isomorphic to a dual affine plane,
then for all $\lambda$ and $\mu\neq 0$ in $\K_{\sigma,-1}$
we have $\lambda\mu^{-1}\in \K_{\sigma,-1}$. 
But that implies that $\K_{\sigma,-1}$ is a 
(commutative) field.
Conversely, if $\K_{\sigma,-1}$ is a commutative field,
then all points of the form $\langle  v_1\lambda+ v_2\mu+v_3\mu\rangle$
with $\lambda,\mu,\rho\in \K_{\sigma,-1}$ different from
$\langle v_1-v_3\rangle$ form a subspace of $\Pi$ isomorphic to
a dual affine plane.

Now consider the case that
$v_1$, $v_2$ and $v_3$ are three isotropic vectors
generating a $3$-space without isotropic lines.

Then without loss of generality
we have $f(v_1,v_2)=1=f(v_1,v_3)$
and $f(v_2,v_3)=\alpha\neq 0$.
Since there is no
isotropic $2$-space in $\langle v_1,v_2,v_3\rangle$,
we find that 
$$f(v_1+v_2,v_1\lambda +v_3)=1-\lambda+\alpha\neq 0$$
for all $\lambda\in \K_{\sigma,-1}$. This 
is equivalent with  $\alpha\not \in \K_{\sigma,-1}=\K^{\sigma,-1}$.

Consider an isotropic point  $\langle v\rangle\subseteq \langle v_1,v_2\rangle$
and an isotropic point $\langle w\rangle\subseteq \langle v_1,v_3\rangle$.
Then there are $\lambda,\mu\in \K_{\sigma,-1}\setminus\{0\}$ with
(up to a scalar) 
$v=v_1+v_2\lambda $ and $w=v_1+ v_3\mu$.
The intersection point
of $\langle v,w\rangle$ with $\langle v_2,v_3\rangle$
is then the point
$\langle  v_2\lambda -v_3\mu \rangle$.
Now we have
$$f( v_2 \lambda-v_3\mu , v_2\lambda - v_3\mu)=
-\lambda^\sigma\alpha\mu+\mu^\sigma\alpha^\sigma\lambda=
(\lambda\alpha\mu)^\sigma-\lambda\alpha\mu.$$
Clearly, if $\K_{\sigma,-1}=\K^{\sigma,-1}$ is a field, 
then  $\lambda\alpha\mu\not \in \K^{\sigma,-1}$ and we find the intersection point not to be isotropic.
This implies that in $\langle v_1,v_2,v_3\rangle$
there are no $4$ hyperbolic lines meeting in $6$ isotropic
points.

Now assume that $\mathcal{S}=\mathcal{S}_q$ for some non-degenerate 
quadratic or
pseudo-quadratic form $q$ with associated
 $(\sigma,-1)$-hermitian
form $f$. 
Suppose $\langle v_1\rangle$ and $\langle v_2\rangle$ are two isotropic points
with $f(v_1,v_2)=1$.
Then the points on the hyperbolic
line through these two points are, besides $\langle v_2\rangle$,
the unique isotropic points in each of the subspaces 
$\langle v_1\lambda+v_2\rangle+\Rad{W,f}$, where $\lambda\in \Phi_q$.
So, hyperbolic lines of $\mathcal{S}$ are thick
if and only if the set  $\Phi_q$ contains at least $2$ elements.
Suppose $\Phi_q$ contains at least $2$ elements.
Let $\mathcal{H}$ be the set of all hyperbolic lines and set $\Pi:=(\mathcal{P},\mathcal{H})$. 
As above we can prove that  any triple of points $p,q,r$ with
$p\sim q\sim r\perp p$ are in a subspace of $\Pi$ isomorphic to
a dual  affine plane if and only if 
for all $\lambda,\mu\in \Phi_q$ we have $\lambda\mu^{-1}\in \Phi_q$.
Since $1\in \Phi_q$, we find $\Phi_q$ to be a commutative field.
Moreover, the same arguments as used above also reveal that,
provided  $\Phi_q$ is a commutative field,
in a subspace spanned by two intersecting
hyperbolic lines but not containing an isotropic line,
there are no $4$ hyperbolic lines meeting in $6$ isotropic points.
}\end{nummer}

\medskip

{}From the above we now easily deduce the following.

\begin{nummer}{\rm\bf Proposition.}\label{classicalplane} 
Let $\mathcal{S}$ be a non-degenerate classical polar space ${\mathcal S}_f$
or $\mathcal{S}_q$ 
defined by a $(\sigma,\epsilon)$-hermitian form $f$
as in {\rm \Ref{classicalpolar}(I)}
or a (pseudo-)quadratic form $q$ with associated $(\sigma,\epsilon)$-hermitian
form $f$ as in {\rm \Ref{classicalpolar}(II)} or {\rm (III)}, respectively.

Suppose $\mathcal S$ contains lines
and let $\Pi$ be the corresponding geometry of hyperbolic lines.
Then $\Pi$ is a connected and co-connected 
thick partial linear space satisfying  the conditions
{\rm (a)-(d)} of Theorem $\ref{planethm}$ if and only if 

\begin{enumerate}
\item $f$ is  as in {\rm \Ref{classicalpolar}(I)} and 
$\K_{\sigma,\epsilon}$ is a commutative field of order at least $3$;
\item  $q$ is 
as in {\rm \Ref{classicalpolar}(II) } or {\rm (III) }
and 
$\Phi_q$ is a commutative field of order at least $3$.
\end{enumerate}
\end{nummer}

\begin{nummer}{\rm {\bf Exceptional polar spaces of rank $3$.}
Besides the classical polar spaces there are many classes of generalized
quadrangles and two more classes of polar spaces of rank $3$; see \cite{tits}.
We now consider the two additional classes of polar spaces of rank $3$.

Let $\mathbb{P}$ be a projective space of (projective) dimension
$3$. Then denote by $P$ the set of lines of $\mathbb{P}$ and by
$L$ the set of pencils of lines on a point inside a plane of $\mathbb{P}$.
Then $(P,L)$ is a polar space. If the space $\mathbb{P}$ is defined over
a field, then by the Klein correspondence we can identify 
this polar space with an orthogonal polar space as described above.
However, if $\mathbb{P}$ is defined over a proper skew field, 
then this provides us with a new polar space.
If $l$ and $m$ are two skew lines of $\mathbb{P}$,
then $\{l,m\}^\perp$ consists of  all the lines of $\mathbb{P}$
meeting both $l$ and $m$.
If $n$ is a line that meets all lines in $\{l,m\}$,
then it is straightforward to check that $n=l$ or $m$.
Thus all hyperbolic lines of the polar space $(P,L)$ are thin.

The second class of polar spaces of rank $3$
are the exceptional polar spaces of type $E_{7}$
over some field $k$ as described in \cite[Section 9]{tits}.
It has been shown in \cite{timmes} that
the automorphism group of this polar space contains a class
of $k$-transvection groups whose elements are in one-to-one
correspondence to the points of the polar space.

As follows from \cite{k-trans}, the geometry of these $k$-transvection
subgroups is a geometry satisfying the conditions of Theorem \ref{main}.
In fact, in this geometry the subspace generated
by three points $x\sim y\sim z\perp x$ is a dual affine plane.
Unfortunately, it is not known to us what the structure of
an arbitrary plane is.
}\end{nummer}

\bigskip

We finish this section with a description of Fischer spaces, some of which
appear in the conclusion of Theorem \ref{fischer}.

\begin{nummer}{\rm\bf Fischer spaces.}\label{fischer-ex}
{\rm 
A {\em Fischer space} is a partial geometry $\Pi=(P,L)$
of order $2$ (i.e., all lines contain $2+1=3$ points)
 in which any pair of intersecting lines is contained in a subspace
isomorphic to a dual affine plane or an affine plane.
Fischer spaces are closely related to groups generated by $3$-transpositions.
Indeed, if $D$ is a normal set of $3$-transpositions in a group $G$, that
is a normal set of involutions with the property that $de$ is of order at most
$3$ for any two elements $d,e\in D$, then the space with $D$ as point set and
as lines the triple of elements from $D$ in the subgroups
of  the form $\langle d,e\rangle$ with $d,e\in D$ with $de$ of order $3$ is a 
Fischer space.
On the other and  if $\Pi$ is a Fischer space, then to each point $p\in P$
we can attach an involutory automorphism $\tau_p$ of $\Pi$ which fixes $p$ and
all points not collinear to $p$ and switches all points $q,r$ with $\{p,q,r\}$
being a line.
The set $\{\tau_p\mid p\in P\}$ is a normal set of $3$-transpositions in the
group
${\rm Aut}(\Pi)$.

By exploiting this relationship all connected Fischer spaces containing
two non-collinear points have been classified by Cuypers and Hall \cite{ch-class,ch-fischer}.
We briefly describe the so-called irreducible Fischer spaces (see \cite{ch-class,ch-fischer}). 

\begin{enumerate}
\item {\em Symmetric Fischer spaces}. Let $\Omega$ be a set and $D$ the class
  of transpositions in $G={\rm Sym}(\Omega)$ the symmetric group on $\Omega$.
Then $D$ is a normal set of $3$-transpositions in $G$. 
The corresponding  Fischer space is called a  symmetric Fischer space.
\item {\em Symplectic and unitary Fischer spaces}. These are the geometries of
  hyperbolic lines of a non-degenerate symplectic polar space over the field
  $\F_2$ or a unitary polar space over $\F_4$.    
\item {\em Orthogonal Fischer spaces over $\F_2$}. Let $q$ be a non-degenerate quadratic form
on a vector space $V$ over $\F_2$ with trivial radical. Then the space whose
point set is the set of all non-isotropic $1$-spaces and whose lines
are the elliptic lines is an orthogonal  Fischer space over $\F_2$. 
\item {\em Orthogonal Fischer spaces  over $\F_3$}. Let $q$ be a non-degenerate quadratic form
on a vector space $V$ over $\F_3$. Take the space with as points the
non-isotropic $1$-spaces and as lines the sets of three points contained
in a 2-space tangent to the quadric. Then  any connected component of this
geometry is called an orthogonal  Fischer space over $\F_3$.
\item {\em Sporadic Fischer spaces}. These are the Fischer spaces
corresponding to 
the classes of $3$-transpositions in a sporadic Fischer group
${\rm Fi}_{22}$, ${\rm Fi}_{23}$, ${\rm Fi}_{24}$, $\Omega(8,2):{\rm Sym}_3$ or 
$\Omega(8,3):{\rm Sym}_3$. 
\end{enumerate}
}\end{nummer}

\section{A Reduction Theorem}
\label{reduction-sec}

In this section we start with a proof of  Theorem \ref{main}. 
Actually we first consider a more
general situation as in the following:

\begin{nummer}{\bf Setting.} \label{setting}
Let $\Pi=(P,L)$ be a  connected partial linear space 
satisfying the following conditions.

\begin{enumerate}
\item All lines in $L$ contain at least $4$ points;
\item if $x$ and $y$ are points with $x^\perp \subseteq y^\perp$ then
$x^\perp=y^\perp$;
\item there exist non-collinear points;
\item every triple of points $x,y,z\in P$ with $x\sim y\sim z$ and $x\perp z$
is contained in a transversal plane;
\item if $\pi$ is a transversal 
plane of $\Pi$ 
and $p$ a point, then $p^\perp \cap \pi$
is empty, a point, a line, a transversal of $\pi$ or $\pi$ itself. 
\end{enumerate}
\end{nummer}

\medskip

It is the purpose of this section to show that we can reduce the above setting 
to the setting of  Theorem \ref{main}.

Let $\Pi$ be as in \Ref{setting}.

\begin{nummer}\label{delta}
If $x$ is a point and $l$ a line of $\Pi$, then $x$ is collinear with 
no point on $l$, with all but one point on $l$ or with all points of $l$.
\end{nummer}

\begin{proof} Suppose $x$ is collinear to some but not all points of the line
  $l$. Then we can find a point $y\in l\cap x^\sim$ and a point $z\in l\cap
  x^\perp$. The points $x,y,z$ form a triangle and hence are contained in
  a transversal plane $\pi$ of $\Pi$.
Clearly $l\subseteq \pi$. Within $X$ we see that the point $x$ is collinear to
all but one of the points of $l$.
\end{proof}

\bigskip

The above shows that $\Pi$ is a so-called $\Delta$-space; see \cite{higman}.

\begin{nummer}\label{transversal}
Let $\pi$ be a transversal plane in $\Pi$.
Then every transversal of $\pi$ is a coclique of size at least $3$.
\end{nummer}

\begin{proof}
See \cite{secants}.
\end{proof}

\begin{nummer}\label{perpissubspace}
Let $x\in P$ be a point then $x^\perp$ and $x^\perp\setminus\{x\}$
are subspaces of $\Pi$.
\end{nummer}

\begin{proof}
Let $l$ be a line in $L$ meeting $x^\perp$ or $x^\perp\setminus\{x\}$ in at least two points. Then 
\Ref{delta} implies that $l$ is contained in $x^\perp\setminus \{x\}$. Thus
$x^\perp$
and $x^\perp\setminus\{x\}$ are  
subspaces of $\Pi$. 
\end{proof}

\begin{nummer}\label{diameter2}
Every connected subspace of $\Pi$ has diameter at most $2$.
\end{nummer}

\begin{proof} Suppose $v,x,y,z$ is a path of length $3$ in the collinearity
  graph of $\Pi$. By \Ref{delta} the points $v$ and $z$  are 
collinear to all but one of the points on the line $l=xy$ on $x$ and $y$.
Since $l$ contains at least $3$ points, there is a point $u\in l$ collinear 
to both $v$ and $z$.
This clearly implies that the diameter of any connected component
is at most $2$.
\end{proof}

\medskip

\begin{nummer}\label{linemeetplane}
Let $\pi$ be a transversal plane and $l$ a line meeting $\pi$ in a point
$x$.
If for some $y\in l\setminus \{x\}$ we have $y^\perp\cap \pi$ being a
transversal   of $\pi$, then for all  
$y\in l\setminus \{x\}$ the intersection  $y^\perp\cap \pi$
is a transversal of $\pi$.

If for some $y\in l\setminus \{x\}$ we have $y^\perp\cap \pi$ being a
line of $\pi$, then for all  
$y\in l\setminus \{x\}$ the intersection  $y^\perp\cap \pi$
is a line of $\pi$.

\end{nummer}

\begin{proof}
Suppose $y$ and $z$ are two points of $l\setminus \{x\}$
with $y^\perp\cap \pi$ being a transversal  $T_y$ of $\pi$.
Let $m$ be a line of $\pi$ on $x$. Then $\langle m,l\rangle$ is a transversal
plane. So $m$ contains a point different from $x$ which is non-collinear to
$z$.
As the above is true for every line $m$ in $\pi$ on $x$, we find that
$z^\perp\cap \pi$ contains at least $2$ points.
But then $z^\perp\cap \pi$ is either a transversal  or a line of $\pi$.
Suppose we are in the latter case and $n=z^\perp\cap \pi$ is a line.
Then this line $n$ meets $T_y$ in some  point, $u$ say.
Inside the transversal plane $\pi$ we see that
all points of  $T_y$, and hence in particular $u$, are  collinear to $x$.
However, $u^\perp$ contains $y$ and $z$ and by \Ref{delta}  all points of $l$
including $x$. This contradiction proves the first statement.

The second follows immediately from \Ref{setting}(e) and the above.  
\end{proof}

\begin{nummer}\label{intersectionofplanes}
Let $x$ and $y$ be two non-collinear
points and suppose $z_1\neq z_2$ are two collinear points both collinear
to $x$ and $y$.

If there is a line $\ell$ in $z_1^\perp$ meeting both $z_2x$ and $z_2y$,
then there is a line $m$ in $z_2^\perp$ meting both $z_1x$ and $z1y$.
\end{nummer}

\begin{proof} 
As $z_1$ is not collinear to all the points on $yz_2$,
we find that $z_1,z_2$ and $y$ generate a transversal plane.
In particular, there is a point $u\in z_1z_2$ which is not collinear to
$y$. 

As $z_1^\perp$ meets the transversal plane
$\langle x,y,z_2\rangle$ in the line $\ell$, the intersection of
$u^\perp$ with this plane is also a line; see \Ref{linemeetplane}.
In particular, $u$ is collinear with $x$.
But, as the line $xz_2$ meets $\ell$, we find that $\langle z_1,z_2,x\rangle$
is a transversal plane, and hence
$u^\perp$ meets $xz_1$ in a point different from $x$. This implies that
$u^\perp$ meets the transversal plane $\langle x,y,z_1$ in a line.
Now \Ref{linemeetplane} implies that also $z_2^\perp$ meets $\langle x,y,z_1$ in a line.
\end{proof}

\begin{nummer}\label{perpdiameter}
The connected components of $(P,\perp)$ have diameter at most $2$.
\end{nummer}

\begin{proof} Let $a\perp b \perp c \perp d$ be a path of length 3 in $(P,\perp)$.
Then $\langle a,c,d \rangle $ is a transversal plane.
Pick a point $x$ on $ac$ different from $a$ and $c$ and
a line on $x$ meeting $ad$
in a point $y$ different from $a$ and $d$. 
This line meets the transversal of
 $\langle a,c,d\rangle$
on $c$ in a point $e$.
Clearly, $b$ is collinear to $e,d$ and $y$.
Now applying \Ref{intersectionofplanes}
yields 
that there is a line $m$ meeting $bd$ and $be$ which is contained
in $y^\perp$.
Let $f$ be the point on $m$ in the transversal 
on $d$ and $e$ in $\langle b,d,e\rangle$.
Then $f^\perp$ contains $dy$  and hence $a$.
 Thus we have found a path of length 2 from $a$ to $d$.
This shows that the diameter of the connected components of
$(P,\perp)$ is at most 2.  \end{proof}

\begin{nummer}{\rm\bf Corollary.}\label{connectedcomps}
Connected components of $(P,\perp)$ of diameter $2$ are connected subspaces
of $\Pi$.
\end{nummer}

\begin{proof} If $x$ and $y$ are two collinear points in a connected component
of $(P,\perp)$, then by the above \Ref{perpdiameter} 
we can find a point $z$ non-collinear to both
$x$ and $y$. But it follows from \Ref{delta} 
that all points on the line $xy$ are in $z^\perp$. Hence $xy$ is
contained in the connected component of $(P, \perp)$ containing $x$, $y$ and $z$.
This shows that this component is a subspace of $\Pi$. 

It remains to prove connectedness in $\Pi$ of a connected component
$C$ of $(P, \perp)$. Thus let $x$ and $y$ be two non-collinear points in
$C$, and 
$z$ a point collinear to both of them. We may assume that $z$ is in a
connected component of $(P,\perp)$ different from $C$.
Let $u$ be a point in the transversal plane generated by $x$, $y$ and $z$, different
from $x$ and $y$, but on the transversal coclique on $x$ and $y$ inside this
plane.
If $x^\perp\neq u^\perp$, then by condition \Ref{setting}(c) there exists a point
$v\in u^\perp$ collinear to $x$. The point $v$ is in $C$ and, by condition
\Ref{setting}(d), also collinear to $y$. 
Similarly  $y^\perp\neq u^\perp$ implies that there is a point
in $C$ collinear to both $x$ and $y$.

If  $x^\perp=y^\perp=u^\perp$, then, as $C$ has diameter $2$, there is point
$v\in C$ collinear to all three points $x,y$ and $u$.
In any case we have found a point in $C$ collinear to both $x$ and $y$. 
\end{proof}

\begin{nummer}\label{onecomponent}
Suppose $(P,\perp)$ contains a connected component of diameter $1$.
Then all components of $(P,\perp)$   have diameter $1$.
\end{nummer}

\begin{proof}
Suppose $C_1,C_2$ are connected components of $(P,\perp)$
with $C_1$ having diameter $1$ and  $C_2$ having diameter $2$. Suppose $x,y,z$ is a path of length 2
in $C_2$ and $u$ is a point of $C_1$. 
Fix a point $p$ on $xu$ different from $x$ and $u$ and
a point $q$ on $uz$ 
different from $z$ and $u$.
The subspaces of $\Pi$ spanned by
$u,y,x$ and by $u,y,z$ are both transversal planes.
Inside these planes we find 
a points $r\in py$ and $s\in qz$ not collinear with $u$.
Clearly $r,s\in C_1$. 
So
$r\perp s$.
The subspace  $\langle y,s,u \rangle$ of $\Pi$  is a
transversal plane. The intersection of
 $r^\perp$ with  $\langle y,s,u \rangle$  equals the  
transversal coclique of  $\langle y,s,u \rangle$ containing $u$ and $s$.
Applying \Ref{linemeetplane} to the line $yr$ and plane 
$ \langle y,q,u \rangle $
implies that  $p^\perp \cap\langle y,q,u \rangle $ is also 
a transversal coclique.
Applying \Ref{linemeetplane}
again but now to  the line $up$ and 
$\langle y,q,u \rangle$ shows that  also $x^\perp$ meets
$\langle u,y,q\rangle$ in a transversal
coclique.
However, as $y$ is in $x^\perp$, 
this  transversal coclique is the unique transversal coclique
 containing  $y$
and  $z$.
This  contradicts that $x$ and $z$ are collinear. \end{proof}

\medskip

On the point set $P$ of $\Pi$ we can define the relation $\equiv$
by $$x\equiv y \Leftrightarrow x^\perp=y^\perp.$$
This relation is clearly an equivalence relation on $P$.
The $\equiv$-equivalence class of a point $x$ is denoted by $[x]$.
For each line $l \in L$ define $[l]$ to be the set $\{ [x]\ |\ x\in L\}$ and 
let $L/\equiv$ 
be $\{ [l]\ | \ l\in L\}$.

\begin{nummer}\label{perponclass}
Suppose $x, y\in P$. Then $x\perp y$ if and only if
$x'\perp y'$ for all $x'\in [x]$ and $y'\in [y]$.
\end{nummer}

\begin{proof}
Suppose $x\perp y$, $x'\in [x]$ and  $y'\in [y]$. 
Then $y\in x^\perp=x'^\perp$.
So $x'\in y^\perp$ which equals $y'^\perp$.
\end{proof}

\medskip
The above Lemma \ref{perponclass}  justifies the following definitions.
If $x,y\in P$, then we write $[x]\perp [y]$ if and only if $x\perp y$.
Moreover,  $[x]^\perp$ 
denotes the set of all classes $[y]$ with $x\perp y$.

\begin{nummer}{\rm \bf Reduction Theorem.}\label{reduction}
$\Pi/\equiv:=(P/\equiv, L/\equiv)$ is a connected partial linear space satisfying:
\begin{enumerate}
\item all lines contain at least 4 points;
\item every triple of points $X,Y,Z$ with $X\sim Y\sim Z\perp X$ is contained in a transversal plane;
\item if $X$ and $Y$ are points with $X^\perp \subseteq Y^\perp$ then $X=Y$;
\item if $\pi$ is a transversal 
plane of $\Pi/\equiv$ 
and $X$ a point, then $X^\perp \cap \pi$
is empty, a point, a line, a transversal of $\pi$ or $\pi$ itself.
\item $(P/\equiv, L/\equiv)$  is linear if and only if $(P,\perp)$ 
has only components of diameter $1$.
\end{enumerate}
\end{nummer}

\begin{proof}
If $p$ and $q$ are collinear points, then $[p]\neq [q]$, so all elements
of $L/\equiv$ contain at least 4 points of $P/\equiv$.
This proves (a) and connectedness of $\Pi/\equiv$.

Let $X$ and $Y$ be points of $P/\equiv$. Then $X\perp Y$ if and only
if $x\perp y$ for some points $x \in X$ and $y\in Y$, if and only if $x\perp y$
for all points $x\in X$ and $y\in Y$. Thus if $X^\perp \subseteq Y^\perp$, then for
each $x\in X$ and $y\in Y$ we have $x^\perp \subseteq y^\perp$, hence $x\equiv y$
and $X=Y$. This proves (c).

Now suppose $X$ and $Y$ are two points in $P/\equiv$ on a line $[l]$, $l \in L$.
Then, for any two points $x\in X$ and $y\in Y$ we have $x\sim y$. 
Fix $x\in X$ and $y\in Y$ and let $z$ be a point on the line through $x$ and
$y$.
We will show that $[z] \in [l]$.
First assume that $x \in l\neq xy$. Then $\langle l,x,y \rangle$ is a transversal plane, and $l$ contains a unique point $z'\in z^\perp$. Suppose $u\perp z'$.
If $u\perp x$, then $u\perp a$ for all $a\in \langle l,y \rangle$. Thus $u\perp z$.
If $u\sim x$, then there is a point on $ux$ non-collinear to the unique point in
$[y] \cap l$ and therefore also to $y$ and all points in the transversal 
of $\langle l, y \rangle$ on $y$. But then \Ref{linemeetplane} implies that $u^\perp \cap
\langle l,y \rangle$ is also a transversal, in particular the transversal
 on $z'$ which contains $z$. Thus $u\perp z$. This shows that $z'^\perp
\subseteq z^\perp$ and hence $z'\in [z]$.

Now suppose $xy$ does not meet $l$. Then let $x'$ be the point on $l$ which is
in $[x]$. Consider the line $x'y$. By the above there is a point $z' \in
x'y \cap [z]$ and a point $z'' \in l \cap [z']$. But then $z''\in [z] \cap l$.
Hence, in any case $[l]$ contains $[z]$.

This proves that $(P/\equiv, L/\equiv)$ is a partial linear space.
Furthermore, it shows that if $X,Y,Z$  are three points
with $X\sim Y\sim Z\perp X$ in $(P/\equiv,L/\equiv)$,
then there are points $x\sim y\sim z\perp x$ in $\Pi$ with $(X,Y,Z)=([x],[y],[z])$.
No two distinct points in $\langle x,y,z \rangle$
are $\equiv $-equivalent. Indeed, suppose $u\perp v$ are two distinct 
points in $\langle x,y,z\rangle$ with $u\equiv v$. Then $u\equiv v'$
for every point $u'$ on the transversal of $\langle x,y,z\rangle$
containing $u$ and $v$. So, without loss of generality, we can assume
$u$ to be on the line through $y$ and $x$ and $v$ to be on the line 
through $y$ and $z$. Moreover, since $[x]\neq [z]$, we find that
$u\neq x$ and $v\neq z$. 
Let  $w\in x^\perp\setminus z^\perp$. 
Then $w\not \in u^\perp$ for otherwise it is in $x^\perp\cap u^\perp\cap
v^\perp$, which is contained in $z^\perp$.
So, $\langle w,y,x\rangle$ is a transversal plane containing
a point  $w'\neq u$ which is in $u^\perp$ and thus also in $v^\perp$.
So $w'^\perp$ meets $\langle x,y,z\rangle$ in a transversal and,
by \Ref{linemeetplane}, so does $w^\perp$, which contradicts
that $w\not \in z^\perp$.
 Hence $\{ [u]\ |\ u\in \langle x,y,z \rangle\}$ is
a subspace of $(P/\equiv, L/\equiv)$ isomorphic to $\langle x,y,z \rangle$ and
hence a transversal plane. Furthermore, (d) follows now easily by
the above.

Finally,  $(P/\equiv, L/\equiv)$ is linear if and only if 
$x^\perp=y^\perp$ for all $x,y\in P$ with $x\perp y$. 
But that is equivalent with $(P,\perp)$ having only components of
diameter $1$. Now \Ref{onecomponent}
implies (e).
 \end{proof}

\section{The geometry of hyperbolic lines}
\label{geometry-sec}

Now we focus on the proof of Theorem \ref{main}.
So, assume that $\Pi=(P,L)$ is a partial linear space satisfying the
hypothesis of that theorem. 

\begin{nummer}\label{perpconnected}
Let $x$ be a point in $P$. The subspace $x^\perp \setminus \{x\}$ of $\Pi$ is
connected.
\end{nummer}

\begin{proof} Suppose $y$ and $z$ are two non-collinear points in $x^\perp \setminus \{x\}$.
Since $\Pi$ is connected and of diameter at most 2,
there is a point $u$ collinear with $y$ and $z$. If $u\perp x$ then we are done,
thus suppose that $u$ and $x$ are collinear. 
Then $\langle x,u,y \rangle$,
$\langle x,u,z \rangle$ and $\langle y,u,z \rangle$ are transversal planes.
By condition \Ref{main}(c)  
and our assumption that all $\equiv$-classes are of size 1, there  is a point $v$ collinear with $y$ but
not with $x$. If $v$ and $z$ are collinear, we are done. Assume $z\perp v$.
Thus by condition \Ref{main}(d) the point $z$ is the only point on the transversal coclique $T$
on $z$ in the transversal
plane $\langle y,u,z \rangle$ not collinear with $v$.

Again by condition \Ref{main}(c), 
there is a point $w$ in $x^\perp$ collinear with
$z$. By the same reasoning as above we can assume that $y$ is the unique point on $T$
not collinear with $w$. Since $T$ contains at least 3 points, there is a point
on $T$ collinear with $v$ and $w$, and by condition \Ref{main}(d) not to $x$. But that means
that we have found a path from $y$ to $z$ in the collinearity graph of $\Pi$
that is completely inside $x^\perp \setminus \{x\}$. \end{proof}

\medskip

The following lemma will be of great importance in constructing a polar space
on the point set $P$.

\begin{nummer}\label{oneorall}
Let $x$ and $y$ be two non-collinear points and suppose $z$ is collinear
to both $x$ and $y$. Then there is a transversal  $T$ on $x$ and $y$ that contains a  point in $z^\perp$.
\end{nummer}

\begin{proof} Let $x$ and $y$ be two non-collinear points and $z$ a point collinear to both of
them. Then by \Ref{perpdiameter} there is a point $u$ in $x^\perp \cap z^\perp$.

First suppose that $u$ and $y$ are collinear. 
Let $u'$ be a point on $yu$ different from $u$ and $y$. 
Then $x$ is not collinear to $u'$.
Moreover, there is a point $z'$ on $yz$ not collinear to $u'$.
This point is different from $y$ and hence collinear with $x$.
Consider now line on $z$ meeting $xz'$ in a point different from $x$ and $z'$.
This line contains a point $w$ in the transversal of $x$ inside the plane
$\langle x,y,z\rangle$.
Notice that $w$ and $u'$ are collinear.
As $u'^\perp$ contains a point of the line $xz'$,
we find $\langle u', w,z\rangle$ to be a transversal plane.
So, there is a point $p$ on the line $zu'$ which is not collinear
to $w$.
Since $x$ is collinear to $z$ and not $u'$, it is collinear to $p$.
By condition \Ref{setting}(e), also  $y$ is collinear to $p$. 
So $w^\perp $ meets
the transversal plane generated by $y$, $z$ and $u$ in the line $yp$.
Fix the point $t$ of this line that is not collinear with $z$.
Then $t$ is collinear with $x$ and $y$ but not with $z$ and $w$.
Now \Ref{intersectionofplanes}
applied to $x,y,t$ and $z'$ implies that there is a line
in $\langle x,y,t\rangle$ which is contained in $z'^\perp$.
This line contains a point $q$ which is in the transversal of $x$ and $y$
inside this plane.
But, then $q^\perp$ contains $y$ and $z'$ and hence also $z$.
So we have found the point we are looking for.

Now assume that $u\perp y$. By the previous 
lemma there is a point $v$ in $x^\perp$
collinear with both $u$ and $y$. If $v\perp z$ we are in the above situation
with $v$ instead of $u$. Thus we can assume that $v$ and $z$ are collinear.
By the arguments in the preceding paragraph, but now with $v$ instead of $y$,
there is a transversal coclique on $x$ and $v$
containing a point $u'$ in $z^\perp$. Since $y$ and $v$ are collinear but $y$ and $x$ are not, we have
that $y$ and $u'$ are collinear. So we can apply the above with $u'$ instead of $u$. 
\end{proof}

\begin{nummer}\label{equivrelation}
Let $x$ be a point and $y$ and $z$ two points in $x^\perp$ different from $x$.
If $x^\perp \cap y^\perp \subseteq x^\perp \cap z^\perp$, then $x^\perp \cap y^\perp=
x^\perp \cap z^\perp$.
\end{nummer}

\begin{proof} 
Suppose $x^\perp \cap y^\perp$ is a proper subset of $x^\perp \cap z^\perp$.
Then $y\perp z$. 

By \Ref{perpconnected} and \Ref{diameter2} there is a point $u\in x^\perp$ collinear to both $y$
and $z$.

Now let $v$ be a point in $x^\perp \cap z^\perp$ collinear with $y$.
By the above Lemma \ref{oneorall}, 
there is a transversal  on $z$ and $v$ containing
a unique point $t$ (possibly equal to $v$)
in $u^\perp$. This point $t$ is in $x^\perp \cap z^\perp$,
but collinear with $y$. Applying \Ref{oneorall} one more time, we find
a transversal  on $u$ and $t$ containing a point $s$ in $y^\perp$.
Then $s$ is in $x^\perp\cap y^\perp$, but $s$ is collinear with $z$, since $z\perp t$ but $z\sim u$.
A contradiction, which proves the lemma. 
\end{proof}

\medskip

Now we are in a position to 
define a set $\cal L$ on $P$ such that $(P, \cal L)$
is a polar space.

Let $x$ and $y$ be two distinct points of $\Pi$. Then define

$$ (xy)^{\perp\perp} =\{ z\in P\ |\ z \perp u \ {\rm for\ all}\ u\in x^\perp \cap y^\perp \}.$$

Now let $\cal L$ be the set of all subsets $(xy)^{\perp\perp}$ of $P$, with $x$ and $y$
two non-collinear points in $P$. The elements of $\cal L$ will be called {\it singular lines}.
If $x\sim y$, then $(xy)^{\perp\perp}$ is called the {\it hyperbolic line} through $x$ and $y$.

\begin{nummer}\label{singularlines}
Let $x$ and $y$ be two non-collinear points. Then for all distinct points $s$ and $t$
in $(xy)^{\perp\perp}$ we have $s\perp t$ and $(st)^{\perp\perp}=(xy)^{\perp\perp}$.
\end{nummer}

\begin{proof} Suppose $s\in (xy)^{\perp\perp}$ different from $x$. Then $s\perp x,y$.
By \Ref{equivrelation} we  have  $x^\perp \cap y^\perp= x^\perp \cap s^\perp$. Hence
$(xy)^{\perp\perp}=(xs)^{\perp\perp}$.
Applying this argument one more time we see $(xs)^{\perp\perp}=(st)^{\perp\perp}$. \end{proof}

\begin{nummer}{\rm \bf Proposition.}\label{ispolar}
$(P,\cal L)$ is a non-degenerate polar space of rank at least $2$.
\end{nummer}

\begin{proof}
By the previous lemma and the definition of singular lines we see that for
any singular line $l$ and point $x$ we have $x^\perp \cap l$ is empty,
consists of a single point or equals $l$. In view of the `one-or-all'
axiom and condition (c) of Theorem \ref{main}, we only have to show that this intersection
is never empty.
So fix two points $y$ and $z$ in $l$ collinear with $x$.
Then $l=(yz)^{\perp\perp}$ and by \Ref{oneorall}
there is a transversal 
on $y$ and $z$ containing  a point in $x^\perp$. By condition (d) of Theorem
\ref{main} this transversal
 is contained in $l$. Thus $x^\perp$ meets $l$ non-trivially.

Since $\Pi$ contains non-collinear points, we find $(P,\mathcal{L})$ to
contain lines. This implies that $(P,\mathcal{L})$ has rank at least $2$. \end{proof}

\medskip

Proposition \ref{ispolar}
finishes the proof of Theorem \ref{main}.

\section{Planar spaces of hyperbolic lines}
\label{plane-sec}

In this section we prove Theorem \ref{planethm}. So, let $\Pi=(P,L)$ be a partial
linear space as in the hypothesis of \Ref{planethm}.

\begin{nummer}
Let $\pi$ be a linear plane of $\Pi$. Any three non-collinear points in $\pi$ generate $\pi$.
\end{nummer}

\begin{proof} Let $x$, $y$ and $z$ be three non-collinear points in $\pi$.
Then the lines $xy$ and $xz$ are distinct and contained in the planes $\pi$
and $\langle x,y,z \rangle$. Since $\Pi$ is planar, these planes are the same. \end{proof}

\medskip

A dual affine plane is transversal.

\begin{nummer}\label{perpmeetsplane}
Let $x$ be a point in $\Pi$ and $\pi$ a dual affine plane.
Then $x^\perp$ meets $\pi$ in the empty set,
a point, a line, the full plane or a transversal.
\end{nummer}

\begin{proof} 
As all planes of $\Pi$ are $\Delta$-spaces, so is $\Pi$.  
In particular, $x^\perp$ is a subspace for any point $x\in P$.
As a dual affine plane  is generated by any three points
not on a line or a transversal, we easily deduce that
for any point $x$ and dual affine plane $\pi$ of $\Pi$ we have that
$x^\perp$ meets $\pi$ in the empty set a point, a line, the full plane,
or part of a transversal.

Suppose that $x^\perp$ meets $\pi$ in at least two points
of a transversal $T$.
Suppose $y$ and $z$ are two
points of $T \cap x^\perp$. If $T$ is not contained in $x^\perp$, then
there is a line $l$ in $\pi$ completely in $x^\sim$. So $x$ and $l$ generate a
linear plane. Fix two points $u$ and $v$ on $l$ but not in $T$.
Let $w$ be the intersection point of $yu$ and $zv$. Let $w'$ be the point
on $l$ not collinear with $w$. Since $x$ is not collinear with $y$, there is a point
$y'$ on $xu$ not collinear to $w$. As $x \perp z$, there is also a point $z'$ on
$xv$ not collinear to $w$. By condition (b) of Theorem \ref{planethm}, the lines $w'y'$ and $w'z'$ are distinct.
Hence $\langle x,l \rangle = \langle w',y',z' \rangle \subseteq w^\perp$, which
is clearly a contradiction. Thus $T \subseteq x^\perp$. \end{proof}

\medskip

With the help of  Theorem \ref{main} and Section \ref{reduction-sec} we obtain the following.

\begin{nummer}{\rm \bf Proposition.}
$(P,\perp)$ is a  non-degenerate polar graph
of rank at  least $2$. 
\end{nummer}

\begin{proof}
By the above results $\Pi$ satisfies the conditions of Theorem \ref{main}.
This proves the result. 
\end{proof}

\medskip

Theorem \ref{planethm} is a direct consequence of the above result.

\section{Hyperbolic lines in polar spaces over finite fields}

\label{finite-sec}

The purpose of this section is to provide a proof of Theorem \ref{fischer}.

We start with the case that $\Pi=(P,L)$ is a partial linear space of order $q\geq 3$
satisfying the hypothesis of Theorem \ref{fischer}.
Since classical unitals are linear spaces in which no $4$ lines meet in $6$ points,
we find that $\Pi$ does satisfy the conditions of
(a)-(c) of Theorem \ref{planethm}.
In particular, we can apply some of the results from the previous sections.

Indeed, by \ref{perpmeetsplane} we can apply the Reduction Theorem
\ref{reduction} and obtain the following.

\begin{nummer}{\bf Proposition}
Suppose $\Pi$ satisfies the conditions of Theorem $\ref{fischer}$.
If $q\geq 3$, then every connected component of
$(P,\perp)$ is a non-degenerate polar graph.
\end{nummer}

We will analyze this situation somewhat further.
Assume that $(P,\perp)$ is connected.
Let $\cal L$ be the set of singular lines.
We first consider the case where
$(P,\cal L)$ is a generalized quadrangle. Since every point in $P$ is on
at least $q+1\geq 4$ singular lines and every singular line contains at least
$q\geq 3$ points, 
all singular lines have the same number of points, and every point
is on the same number (possibly infinite) of singular lines.
Suppose the generalized quadrangle $(P,\cal L)$
has order $(s,t)$, where $s,t\in\N\cup \{\infty\}$.

\begin{nummer}\label{sandtfinite}
$s$ and $t$ are finite.
\end{nummer}

\begin{proof} First we consider $t$.
Suppose there are more than $q+1$ singular lines on a point $x$. Then
$x^\perp \setminus \{x\}$ contains a linear plane $\pi$ say. Let $y$ be a point
on a singular line
through $x$, then $y^\perp \cap \pi$  contains a point, $z$ say, and $y\perp x,z$.
Hence $y$ is contained in the singular line through $x$ and $z$.
This implies that there are at most $|\pi|$ singular lines through $x$ and $t$
is finite.

Now consider $s$.
Suppose there are more than $q+1$ points on a singular line $l$.
Fix distinct points $x$, $y$  and $z$ on $l$ and let $m$ be an ordinary 
line through
$y$ contained in $x^\perp$. Let $n$ be a singular line on $z$ different
from $l$. Then $n$ does not meet $m$ and contains a point $x'$ that is collinear to
all points on $m$. Let $\pi$ be the linear plane generated by $x'$ and $m$.
For every point $p$ of $\pi\setminus \{y\}$ 
there is a unique point in $p^\perp\cap l$.
Thus either $l$ contains at most $|\pi|$
points and hence $s$ is finite, or it 
contains a point $y'$ different from $y$ that is collinear with all points
of $\pi$ except for $y$. Suppose we are in the latter situation.
Then let $k$ be a singular line on $y'$ different from $l$. Then again,
for every point $p\in \pi$,
there is  a unique point in $k\cap p^\perp$.
By assumption \Ref{fischer}(d),  all points of $k$
are collinear (in $\Pi$) to some point of $\pi$. So, $|\pi| \geq |k|$ and $s$ is finite. \end{proof}

\bigskip

Next we show that the elements of $L$ are the hyperbolic lines of the
generalized quadrangle $(P, \cal L)$.

\begin{nummer}
Let $x$ and $y$ be two collinear points and $X$ a set contained in $(xy)^{\perp\perp}$.
Then $\langle X\rangle$ is a linear subspace contained in $(xy)^{\perp\perp}$.
\end{nummer}

\begin{proof} This is a direct consequence of the Delta space property of $\Pi$,
and the fact that in a generalized quadrangle hyperbolic and singular lines
meet in at most one point. \end{proof}

\begin{nummer}
The elements of $L$ are the hyperbolic lines of $(P,\cal L)$.
\end{nummer}

\begin{proof} Let $x$ and $y$ be two collinear points, then by the above lemma
$xy$ is contained in the hyperbolic line $(xy)^{\perp\perp}$. Now suppose there
is a point $z$ in $(xy)^{\perp\perp}$ which is not in $xy$.
Then the plane $\pi=\langle x,y,z \rangle$ is contained in the hyperbolic line
on $x$ and $y$. Let $p$ and $p'$ be  distinct points in $x^\perp \cap
y^\perp$. 
If $l$ is a singular line on $p$ (or $p'$), then
for any point $q\neq p$ on $l$ but not in $\pi$, 
we find, by condition (d) of Theorem \ref{fischer}, a point $r\in q^\perp\cap
\pi\subseteq p^\perp$. So, $r\in l$.
So, all singular lines on $p$ (and also on $p'$) 
meet $\pi$, which therefore has cardinality $t+1$.
By 1.4.2 of \cite{pt} we have that  $s \geq t$. Thus on the singular line $(zp')^{\perp\perp}$
we can find a point $z'$
not in $\pi$ but collinear with all points on $xy$.
Consider the plane $\pi'=\langle x,y,z' \rangle$. This is contained in
$p'^\perp$ and contains therefore at most $t+1$ points. Now every point of $\pi'$
is non collinear with a unique point on $(pz)^{\perp\perp}$
and every point of $(pz)^{\perp\perp}$ is
non collinear to some point in $\pi'$. Since $xy$ is contained in $p^\perp$
this implies that $(pz)^{\perp\perp}$ contains at most $|\pi'|-q$ points.
Hence $s+1\leq |\pi'|-q<t$, a contradiction. \end{proof}

\begin{nummer}{\rm \bf Proposition.}
Suppose $(P,\perp)$ is connected. If  $(P,\cal L)$ has rank $2$, then it is a classical symplectic
or unitary generalized quadrangle defined over the field $\F_q$
or  $\F_{q^2}$, respectively.
\end{nummer}

\begin{proof} By the above $(P, \cal L)$ is a finite generalized quadrangle in which two intersecting
hyperbolic lines generate a linear subspace of $\Pi$ or a dual affine plane. Now Theorem 5.6.8
of \cite{pt} proves the proposition. \end{proof}

\medskip

Now suppose $(P,\cal L)$ is a polar space of rank at least 3 and possibly
infinite.

\begin{nummer}{\rm \bf Proposition.}\label{quad}
Suppose $(P,\perp)$ is connected. Then $(P,\cal L)$ is a classical symplectic
or unitary polar space defined over the field $\F_q$
or  $\F_{q^2}$, respectively.
\end{nummer}

\begin{proof}
If the rank $r$ of $(P,\cal L)$ is finite, then by considering the subspaces $x^\perp \cap y^\perp$ for two points $x$ and $y$ with $x\sim y$,
we obtain a polar space of rank $r-1$, which by the Delta space property of $\Pi$
is also a subspace of $\Pi$. Thus by induction
all singular lines contain
either $q+1$ points, or $q^2+1$ points. Now it follows easily from
the classification of polar spaces of rank at least 3, see \cite{johnson,tits}, that $\Pi$ is one of the spaces
of the conclusion of Theorem \ref{fischer}.

If the rank of the polar space is infinite, then it is embeddable in a projective
space such that hyperbolic lines are contained in projective lines, see \cite{cjp,johnson,tits}. Assume
we have such an embedding of $(P, \cal L)$.
By intersecting the polar space with an appropriate finite dimensional
projective subspace we can get a non-degenerate polar subspace of finite rank which is also a subspace
of $\Pi$. So again we find that singular lines are of size $q$ or $q^2$ and
$L$ is the set of hyperbolic lines in $(P, \cal L)$ and it
follows from
the classification of polar spaces of rank at least 3, (see \cite{johnson,tits}), that $\Pi$
is isomorphic to one of the spaces in the conclusion of Theorem \ref{fischer}.
\end{proof}

\medskip

To finish the proof of Theorem \ref{fischer} in the case that $q\geq 3$ we still have to prove that
the graph $(P,\perp)$ is connected.
But this is done in the next lemma.

\begin{nummer}
The graph $(P,\perp)$ is connected.
\end{nummer}

\begin{proof}
Suppose $C$ is a connected component of $(P,\perp)$ and $x$ a point not in
$C$.
By the above results, $C$ is the polar graph of a classical symplectic or
unitary polar space. In particular singular
lines of $C$ contain $q+1$ or $q^2+1$ points.

Fix a singular line $l$ of $C$ and a line $m$  of 
$\Pi$ on $x$ meeting $l$ in a point.
If $n$ is a second line on $x$ meeting $l$, then
$\langle m,n\rangle$ is a dual affine plane meeting
$l$ in a transversal coclique. As each dual affine plane is generated by any two
of its lines, we find that there are $(|l|-1)/(q-1)$ dual affine
planes on $m$ meeting $l$ in a transversal.
But, that contradicts $q>2$ and we can conclude that $(P,\perp)$ is connected.
\end{proof}

\medskip

This completes the proof of  Theorem \ref{fischer} in the case that $q\geq 3$.
It remains to consider the case that $q=2$.

Suppose that 
$\Pi=(P,L)$ is a partial linear space satisfying the hypothesis
of Theorem \ref{fischer} and assume that all lines of 
$L$ have exactly $3$ points
In that case $\Pi$ is a {\em Fischer space}  and we can apply the results
of Cuypers and Hall \cite{ch-class} and of \cite{chps}.
The conditions in \Ref{fischer}(c)
imply that $\Pi$ is both $\theta$- and $\tau$-reduced in the terminology of 
\cite{ch-class}. Now we easily deduce from
\cite[Section 4]{chps} that $\Pi$ is one of the spaces
in the conclusion of Theorem \ref{fischer}.

\section{Applications}
\label{applications}

In this final section we briefly discuss some applications of our results.
The first is concerned with groups generated by so-called
$\K$-transvections, $\K$ being a field, the second with Lie algebras generated by 
special elements called extremal, 
 and the final one on local recognition
of graphs.

\begin{nummer}{\bf $\K$-Transvection groups.}
{\rm Suppose $\K$ is a field. Let $G$ be a group generated by a
conjugacy class of Abelian subgroups $\Sigma$
such that for any two elements $A,B\in \Sigma$ we have
one of the following:
\begin{enumerate}
\item $\langle A,B\rangle$ is isomorphic to ${\rm (P)SL}(2,\K)$
and $A$ and $B$ are full unipotent subgroups of $\langle A,B\rangle$;
\item $[A,B]=1$, notation $A\perp B$.
\end{enumerate}

Then $\Sigma$ is called a {\em class of $\K$-transvection (sub)groups} of $G$.
Groups generated by $\K$-transvection groups have been studied by
Aschbacher \cite{asch1}, Aschbacher and Hall \cite{ah},
Fischer \cite{fischer,fischer-inv}, Steinbach \cite{invent,steinbach,stein} and others, and in particular 
by  Timmesfeld \cite{timmes,timmes-book}.

Our main result can be viewed as a geometric version of 
some of the results by Timmesfeld. This connection has  been explored  in
\cite{k-trans}.

Let $G$ be a group generated by a conjugacy class of $\K$-transvection groups 
$\Sigma$. We construct a geometry $\Pi(\Sigma)$ 
with point set $\Sigma$. Lines are defined
as follows.
For each non-commuting pair $A,B\in \Sigma$ the line
on $A$ and $B$ is the set of elements $C$ from $\Sigma$ with
$C\leq \langle A,B\rangle$. 

Let $\Sigma$ be a class of $\K$-transvection subgroups of the group $G$.
If $\K$ contains at least $4$ elements, then,
as follows from \cite{k-trans},
we are able to show that  $\Pi(\Sigma)$ satisfies the
following conditions as given in \Ref{setting},
which  implies that we can invoke our main results, Theorem
\ref{reduction} and Theorem \ref{main}. This leads to the following result
(for details, see \cite{k-trans}).
}
\end{nummer}

\begin{nummer}{\bf Theorem.}
Let $\K$ be a field and $G$ be a group generated by a class $\Sigma$ of $\K$-transvection subgroups.
Suppose that $\K$ contains at least $4$ elements and that $G$ has no nontrivial nilpotent 
normal subgroup. If $\Sigma $ contains two commuting elements $A$ and
$B$ with $C_\Sigma (A) \neq C_\Sigma (B)$, then  $(\Sigma , \perp)$
is the polar graph of a non-degenerate polar space of rank at least $2$.
\end{nummer}
 
The above result together with the work of Steinbach \cite{stein,steinbach} and
Timmesfeld
\cite{timmes} provides a full classification of the groups $G$ as
in the hypothesis of the above theorem.

When we allow the field $\K$ to have order $2$ or $3$, the situation
is different. In that case we obtain a geometry in which \ref{main} (d)
might be violated. 
The additional geometries appearing 
in these cases are the so-called {\em generalized Fischer spaces}, see
\cite{c-fischer,k-trans}.
Combining the results of this paper, \cite{k-trans}
and various results on $3$-transposition
groups leads to a classification of all 
groups generated by a class $\Sigma$ of
$\K$-transvection groups, $\K$ an arbitrary field, 
having no nontrivial nilpotent normal subgroups.
For details the reader is referred to \cite{k-trans}.
  
\begin{nummer}{\bf Extremal elements in Lie algebras.}
{\rm  Another situation in which one can apply our main result Theorem \ref{main}
is the following.
Let $\mathcal{L}$ be a Lie algebra over the field $\K$.
An element $x\in \mathcal{L}$ is called {\em extremal}
if $[x,[x,\mathcal{L}]]\subseteq \K x$.
Extremal elements in Lie algebras are studied by Cohen et al. \cite{csuw}
and Cohen and Ivanyos \cite{ci}.
In \cite{csuw} it has been shown that 
for any two distinct extremal elements $x,y\in \mathcal{L}$ we have one of the following:
\begin{enumerate}
\item $[x,y]=0$;
\item $[x,y]$ is an extremal element;
\item $\langle x,y\rangle$ is isomorphic to ${\mathfrak sl}(2,\K)$. 
\end{enumerate}

If case (b) does not occur, then the geometry with as points the $1$-dimensional subspace $\K x$ of
$\mathcal{L}$  
and as lines the sets of points contained in
an ${\mathfrak sl}(2,\K)$-subalgebra generated by two extremal elements
is again a geometry satisfying 
the conditions of the Setting \ref{setting}; see \cite{panhuis}.
Theorem \ref{reduction} and Theorem \ref{main} can now be invoked to obtain
an important step in the classification of Lie algebras generated by extremal
elements. 
 }\end{nummer}

\begin{nummer}{\bf Local recognition of graphs.}
{\rm Our final application is concerned with local recognition of graphs.
Let $\Gamma$ be a graph and $x$ a vertex of $\Gamma$.
By $\Gamma_x$ we denote the subgraph of $\Gamma$ induced on the set of
neighbors of $x$.

Consider the graph $\Gamma=\Gamma(\mathcal{S})$ with as vertices the set  the   hyperbolic lines
of a non-degenerate classical  polar space $\mathcal{S}$, two vertices being adjacent if they
are perpendicular to each other.
Then for each vertex $x$ of $\Gamma$ the induced subgraph  $\Gamma_x$ is isomorphic to $\Gamma(\mathcal{T})$,
where    $\mathcal{T}$ is the polar subspace $x^\perp$ of $\mathcal{S}$.

For example, if $\mathcal{S}$
 is a non-degenerate symplectic polar space of rank $n$, then $\mathcal{T}$
is a symplectic space of rank $n-1$.
It has been proved by
Gramlich \cite{gramlich} that in case $n>5$ the 
graph $\Gamma$ is the unique  connected  graph such that for each
vertex $x$ the graph $\Gamma_x$ is isomorphic to $\Gamma(\mathcal{T})$. 
In the proof of that result, Gramlich uses the main result of 
\cite{c-symp} to identify
$\Gamma(\mathcal{S})$. The Theorems \ref{main} and \ref{planethm} provide tools
to extend Gramlich's methods to obtain similar local recognition
results   on the  graphs of hyperbolic lines of arbitrary
polar spaces. This approach has been applied successfully in \cite{altmann}. 

}\end{nummer}


\bigskip

\begin{thebibliography}{10}

\bibitem{altmann} C. Altmann and R. Gramlich,
On the hyperbolic unitary geometry,
{\it J. Algebraic Combin.} {\bf 31} (2010), 547-583. 


\bibitem{asch1} M. Aschbacher, A characterization of the unitary and symplectic groups
over fields
of characteristic at least 5, {\it Pac. J. Math.} {\bf 47} (1973), 5-26.

\bibitem{ah} M. Aschbacher and M. Hall Jr., Groups generated by a class of elements of order 3,
{\it J. of Algebra} {\bf 24} (1973), 591-612.

\bibitem{buek-Fischer} F. Buekenhout, 
{\it La g\'eom\'etrie des groupes de Fischer}, 
unpublished notes, Free University of Brussels, 1974.

\bibitem{bu-sh} F. Buekenhout and E. Shult, On the foundations of polar
  geometry.  {\it Geom. Dedicata} {\bf   3}  (1974), 155-170. 

\bibitem{csuw}  A.M. Cohen, A. Steinbach, R. Ushirobira, D. Wales, 
{ Lie algebras generated by extremal elements}, 
{\it J. Algebra} {\bf  236}  (2001),  no. 1, 122-154. 

\bibitem{ci}  A.M. Cohen, G. Ivanyos, 
Root filtration spaces from Lie algebras 
and abstract root groups. {\it J. Algebra} {\bf 300} (2006), no. 2, 433-454.

\bibitem{c-fischer} H. Cuypers, On a  generalization of Fischer spaces, {\it Geom. Dedicata} {\bf 34} (1990), 67-87.

\bibitem{c-symp} H. Cuypers, Symplectic geometries, transvection groups and
  modules, {\it J. of Comb. Th.} A, {\bf 65} (1994), 39-59.

\bibitem{c-gen} H. Cuypers, Generalized Fischer spaces, in {\it Finite
  geometry and combinatorics} (Deinze, 1992),  121-129, London
  Math. Soc. Lecture Note Ser., {\bf 191}, 
Cambridge Univ. Press, Cambridge, 1993.

\bibitem{ch-fischer} H. Cuypers and J.I. Hall, The $3$-transposition groups
  with trivial center, {\it J. Algebra} {\bf   178 } (1995),  no. 1, 149-193. 

\bibitem{ch-class} H. Cuypers and J.I. Hall, The classification of 3-transposition groups with trivial center,
in {\it Groups, Combinatorics and Geometry}, proceedings of the L.M.S
symposium on Groups and Combinatorics, Durham, 1990, eds
M. Liebeck and J. Saxl, L.M.S. Lecture Notes Ser. {\bf 165}, 121-138,
Cambr. Univ. Press, 1992.

\bibitem{cjp} H. Cuypers, P. Johnson and A. Pasini, On the embeddability of
polar spaces, {\it Geom. Dedicata} {\bf 44} (1992), 349-358.


\bibitem{cs} H. Cuypers and E.E. Shult, On the classification of generalized Fischer spaces,
{\it Geom. Dedicata} {\bf 34} (1990), 89-99.	

\bibitem{invent} H. Cuypers, A. Steinbach, 
Linear transvection groups and embedded polar spaces,
{\it  Invent. Math.} {\bf   137}  (1999),  no. 1, 169-198. 

\bibitem{k-trans} H. Cuypers, The geometry of $k$-transvections groups, 
{\it  J. Algebra} {\bf 300} (2006), no. 2, 455-471.

\bibitem{secants} H. Cuypers, The geometry of secants in embedded polar spaces, {\it European J. Combin.} {\bf 28} (2007), no. 5, 1455-1472.

\bibitem{chps}
H. Cuypers, M. Horn, J. in 't panhuis and
S. Shpectorov, {\it
Lie algebras and $3$-transpositions}, submitted.

\bibitem{fischer} B. Fischer, {\it Finite groups generated by $3$-transpositions},
University of Warwick lecture notes, 1969.

\bibitem{fischer-inv} B. Fischer, Finite groups generated by 3-transpositions
  I, {\it Inv. Math.} {\bf 13} (1971), 232-246.

\bibitem{gramlich} R. Gramlich,  On the hyperbolic symplectic geometry,
{\it  J. Combin. Theory Ser. A}  {\bf 105}  (2004),  no. 1, 97-110. 


\bibitem{copolar}
J.I. Hall,  The hyperbolic lines of finite symplectic spaces,
{\it  J. Combin. Theory Ser.} A {\bf 47} (1988), no. 2, 284-298.

\bibitem{higman}
D. Higman, Admissible graphs, in {\it Finite geometries}
Proceedings of a Conference in honor of T. G. Ostrom held at Washington State University, Pullman, Wash., April 7-11, 1981. Edited by Norman L. Johnson, Michael J. Kallaher and Calvin T. Long, Pullman, Wash. (1981), pp. 211-222, 
Lecture Notes in Pure and Appl. Math., 82, Dekker, New York, 1983. 

\bibitem{panhuis} J.\ in\ 't\ panhuis, {\it
 Lie algebras, extremal elements, and geometries}, 
Thesis, Eindhoven University of Technology, 2009.\\
{\tt http://alexandria.tue.nl/extra2/200612932.pdf}

\bibitem{johnson} P. Johnson, Polar spaces of arbitrary rank, {\it Geom. Dedicata} {\bf 35} (1990), 229-250.

\bibitem{pt} S. Payne and J. Thas, {\it Finite generalized quadrangles}, Pitman Publ., 1984.

\bibitem{onan} M. O'Nan, Automorphisms of unitary block designs,
{\it J. Algebra} {\bf  20} (1972), 495-511.

\bibitem{steinbach} A. Steinbach, Gruppen, die von  $k$-Transvectionen erzeugt werden,
{\it Mitt. Math. Sem. Univ. Giessen} {\bf 211} (1992), 27-48.

\bibitem{stein} A. Steinbach,  Generalized quadrangles arising from groups
  generated by abstract transvection groups in {\it  
Groups and geometries} (Siena, 1996),  189-199, Trends Math., 
Birkh\"auser, Basel, 1998. 

\bibitem{timmes} F.G. Timmesfeld, Groups generated by $k$-transvections, {\it Inv. Math.}
{\bf 100} (1990), 169-206.

\bibitem{timmes-book} F.G.  Timmesfeld, 
{\it Abstract root subgroups and simple groups of Lie-type},
Monographs Math. {\bf 95}. Birkh\"auser, Basel 2001.

\bibitem{tits} J. Tits, 
{{\it Buildings of spherical type and finite $BN$-pairs}}, 
Lecture Notes in Math. {\bf 386}, Springer Verlag, 1974.

\bibitem{veldkamp} F. Veldkamp,  
Polar geometry, I, II, III, IV, V. {{\it Nederl. Akad. Wetensch. Proc. Ser. A}} {\bf 62} and {\bf  63}  (1959), 512-551, 207-212.
\end{thebibliography}
\end{document}